\documentclass[11pt,twoside,reqno]{amsart}
\usepackage{amsxtra}
\usepackage{color}
\usepackage{amsopn}
\usepackage{amsmath,amsthm,amssymb,arydshln}
\usepackage{mathrsfs,mathtools}
\usepackage{stmaryrd}
\usepackage{hyperref}
\usepackage{enumerate}
\usepackage{xcolor}
\usepackage{float}
\usepackage{pifont}
\usepackage{adjustbox}

\newtheorem*{theorem*}{Theorem}

\theoremstyle{definition}

\theoremstyle{remark}

\makeatletter
\@namedef{subjclassname@2020}{%
  \textup{2020} Mathematics Subject Classification}
\makeatother

\newcommand{\R}{{\mathbb R}}
\renewcommand{\H}{\mathrm{H}}

\newcommand{\beq}{\begin{equation}}
\newcommand{\eeq}{\end{equation}}

\newcommand{\f}{\varphi}

\newcommand{\GL}{{\mathrm {GL}}}

\newcommand{\G}{{\mathrm G}}

\newcommand{\Ng}{{\mathrm N}}

\newcommand{\W}{\wedge}

\newcommand{\frn}{\mathfrak{n}}

\numberwithin{equation}{section}

\title[A homogeneous G$_2$-Laplacian soliton of gradient type]{On the existence of homogeneous solitons of gradient type 
for the G$_{\mathbf2}$-Laplacian flow}

\author{Anna Fino} 
\address{{\scriptsize Dipartimento di Matematica ``G.~Peano'' \\ Universit\`a degli Studi di Torino\\
Via Carlo Alberto 10\\
10123 Torino\\ Italy\\
 and Department of Mathematics and Statistics\\
 Florida International University\\
 Miami  FL 33199, USA}}
\email{annamaria.fino@unito.it, afino@fiu.edu}

\author{Alberto Raffero}
\address{{\scriptsize Dipartimento di Matematica ``G.~Peano'' \\ Universit\`a degli Studi di Torino\\
Via Carlo Alberto 10\\
10123 Torino\\ Italy}}
\email{alberto.raffero@unito.it}


\begin{document}
\begin{abstract} 
In this note, we prove the existence of homogeneous gradient solitons for the G$_2$-Laplacian flow by providing 
the first known example of this type.  
This result singles out the G$_2$-Laplacian flow as the first known geometric flow 
admitting homogeneous gradient solitons on spaces that are one-dimensional extensions in the sense of \cite{PeWi}. 
\end{abstract}
\maketitle

\section{Introduction}

A $\G_2$-structure on a 7-manifold $M$ is given by a positive $3$-form $\varphi$.  
The 3-form $\varphi$ induces a Riemannian metric $g_{\varphi}$ and an orientation on $M,$ and thus  a Hodge star operator $\star _\varphi$. 
Moreover, at each point $x$ of $M,$ there exists a $g_\f$-orthonormal basis $(e^i)_{1\leq i \leq7}$ of $T^*_xM$ such that 
\[
\left.\f\right|_x = e^{127} + e^{347} + e^{567} + e^{135} - e^{146} - e^{236} - e^{245}, 
\]
where $e^{ijk}$ is a shortening for the wedge product of covectors $e^i\W e^j\W e^k$. 
If the $\G_2$-structure is {\em torsion-free}, i.e., $d \varphi = 0$  and $d \star_{\varphi} \varphi = 0$, the metric $g_{\varphi}$ is  Ricci-flat 
with holonomy contained in $\G_2$. 
A $\G_2$-structure satisfying the less restrictive condition $d \varphi =0$ is called {\em closed}. 
In this case, the intrinsic torsion of the structure is encoded into a unique 2-form $\tau$ such that 
\[
d \star_\f \varphi = \tau \wedge \varphi = -\star_\f\tau,  \qquad \Delta_{\varphi}  \varphi = d \tau,
\] 
where $\Delta_{\varphi} = dd^*+d^*d$ is the Hodge Laplacian induced by $\varphi$. 
Closed $\G_2$-structures appear in the study of Riemannian manifolds with holonomy $\G_2$,  
since they are natural candidates to deform toward a torsion-free one.
\smallskip

A closed G$_2$-structure is said to be a {\em Laplacian soliton} if it satisfies the equation
\begin{equation}\label{LS}
\Delta_\f\f = \lambda\f+\mathcal{L}_X\f,
\end{equation}
for some real number $\lambda$ and some vector field $X$ on $M.$  
A Laplacian soliton is said to be of {\em gradient type} if $X=\nabla f = (df)^\sharp$, for some smooth function $f\in C^\infty(M)$.  
Here, the musical isomorphism $(\cdot)^\sharp:T^*M\to TM$ is taken with respect to the metric $g_\f$.
Depending on the sign of $\lambda$, a Laplacian soliton is called  {\em expanding} ($\lambda>0$), {\em steady} ($\lambda=0$),  
or {\em shrinking} ($\lambda<0$). 
On a compact manifold, every Laplacian soliton which is not torsion-free must satisfy \eqref{LS} with $\lambda>0$ 
and $\mathcal{L}_{X}\f\neq0$ (see \cite{Lin,LoWe}).  
The existence of non-trivial Laplacian solitons on compact manifolds is still an open problem. 

Laplacian solitons give rise to self-similar solutions to the Laplacian flow, 
a  geometric flow for closed $\G_2$-structures introduced by R. Bryant in \cite{Bry} and defined by the
equation $\partial_t  \varphi(t) =  \Delta_{\varphi(t)}  \varphi(t)$. 
The short-time existence and uniqueness of the flow was proved in \cite{BrXu} 
(see also \cite{BeVe} for a concise proof based on a result by R.~Hamilton), 
and further geometric and analytic properties of the flow have been studied  in the literature. 
We refer the reader to \cite{Lot} and the references therein for an account on the main results.
\smallskip

In the homogeneous case, i.e., when the automorphism group of $(M, \varphi)$ or a subgroup thereof acts transitively on $M$,   
Laplacian solitons that are not torsion-free may only occur on non-compact manifolds (see  \cite{PoRa}). 
Non-compact expanding, steady and shrinking homogeneous examples have been constructed 
in \cite{FFM,FR1,FR2,Lau1,Lau2,LaNi,LaNi2,Nic1,Nic2}. 
They are given by simply connected Lie groups $\G$ endowed with a left-invariant closed $\G_2$-structure solving \eqref{LS} with respect 
to a vector field $X$ that is not left invariant and is induced by a derivation of the Lie algebra of $\G$. 
Inhomogeneous steady and shrinking Laplacian solitons of gradient type have been provided in \cite{Ball,Fow}. 
More recently, continuous families of complete cohomogeneity-one gradient steady solitons and complete examples of 
cohomogeneity-one gradient shrinking solitons have been constructed in  \cite{HaNo}. 

At the present time, it is not known whether homogeneous Laplacian solitons of gradient type exist. 
In this note, we give an affirmative answer to this problem by showing the following. 
\begin{theorem*}
There exists a simply connected almost nilpotent Lie group diffeomorphic to $\R^7$ 
admitting a left-invariant steady Laplacian soliton of gradient type.  
\end{theorem*}

\smallskip

In the recent work \cite{Ng}, a structure theorem for homogeneous gradient Laplacian solitons has been proved 
using the general structure theorem \cite[Theorem 3.6]{PeWi}. 
The possible situations that may occur  depend on the divergence of the symmetric 2-covariant tensor 
$t_\f\coloneqq g_\f(T_\f^2\cdot,\cdot)$, where 
$T_\f$ denotes the $g_\f$-skew-symmetric endomorphism corresponding to the torsion form $\tau$. 
In detail, assume that $M$ is a 7-dimensional homogeneous space 
admitting an invariant Laplacian soliton of gradient type $\f$, with $X= \nabla f$ and non-constant $f$. 
Then,  $t_\f$ is divergence-free if and only if $(M, g_\varphi)$ is isometric to a product $N \times \R^k$, where $f$ is constant on $N$.  
If $t_\f$ is not divergence-free, then two possibilities occur. 
In the first case, $(M,g_{\varphi})$ is a one-dimensional extension, i.e., $M$ is acted on transitively by a semidirect product $\G\rtimes \R$ 
and is diffeomorphic to $N\times\R$, where $N$ is $\G$-homogeneous, the metric $g_\f$ is of the form 
$g_{\varphi} = g_r +  dr^2$, where $r:M\to\R$ is the distance function to $N$ and $g_r$ is a one-parameter family of invariant metrics on $N,$  
and $f = ar + b$, for some $a,b\in\R$. 
In the second case, $(M, g_{\varphi})$  is isometric to a product $N \times \R^k$, where $N$ is a one-dimensional extension, and 
$f(x,y) = ar(x) + v(y),$ where $v$ is a function on $\R^k$ and $r$ is a distance function on $N$.

\smallskip

A closed $\G_2$-structure  $\varphi$ is said to be {\em extremally Ricci pinched} if its torsion form $\tau$ satisfies the equation 
$d \tau  = \tfrac{1}{6} \vert \tau \vert^2 \varphi + \tfrac{1}{6} \star_{\varphi} \tau^2.$
This condition was introduced in \cite{Bry}, and it characterizes a pinching condition for the Ricci curvature of $g_\f$ on a compact manifold. 
A classification result for homogeneous extremally Ricci pinched closed $\G_2$-structures has been established in \cite{Ball}: 
any such structure is equivalent to a left-invariant one on a simply connected solvable Lie group and it is (locally) a steady Laplacian soliton. 
The classification of extremally Ricci pinched $\G_2$-structures on solvable Lie groups was achieved in \cite{LaNi2}. 

\smallskip

In our previous work \cite{FR2}, we constructed an example of a left-invariant steady Laplacian soliton on a solvable (almost-nilpotent) Lie group $\H$   
satisfying two remarkable properties.
First, it is not an extremally Ricci pinched closed G$_2$-structure. 
Moreover, it satisfies the soliton equation \eqref{LS} with respect to a left-invariant vector field on the Lie group $\H$. 
This last property was not investigated further.  
Here, we show that this homogeneous example is of gradient type by considering an explicit realization of the Lie group $\H$ 
as a matrix subgroup of $\GL(7,\R)$. This proves the main result of the present note. 
Moreover, we prove that this example is a one-dimensional extension, namely it is an instance of 
case 2(a) of the structure theorem \cite[Theorem 1.1]{Ng} recalled above.  
Remarkably, this singles out the G$_2$-Laplacian flow as the first known geometric flow 
admitting examples of gradient solitons on one-dimensional extensions (see \cite[Sect.~1]{PeWi}).

\section{The homogeneous steady Laplacian soliton of gradient type}

Let us consider the subgroup $\H\subset \GL(7,\R)$ whose elements are matrices of the form 
\[
h = \left(\begin{array}{ccccccc}
1 & 0 & 0 & 0 & 0 & 0 & 0 
\\
 0 & 1 & 0 & 0 & 0 & 0 & 0 
\\
 0 & 0 & \exp(x_7) & 0 & 0 & 0 & {x_2}  
\\
 0 & 0 & 0 & \exp(-x_7) & 0 & 0 & -{x_1}  
\\
 2 {x_1}  & 0 & 0 & {x_6}  & \exp(-x_7) & 0 & {x_4}  
\\
 0 & -2 {x_1}  & 0 & {x_5}  & 0 & \exp(-x_7) & {x_3}  
\\
 0 & 0 & 0 & 0 & 0 & 0 & 1 
\end{array}\right),
\]
where $x_i\in\R$, for $1\leq i\leq7$. 
$\H$ is a simply connected, solvable, non-unimodular Lie group of dimension seven. 
$\H$ is diffeomorphic to $\R^7$, and it is isomorphic to a semidirect product of the form $(\R\times \Ng_{5,2}) \rtimes \R$, 
where $\Ng_{5,2}$ is the five-dimensional simply connected 2-step nilpotent Lie group with Lie algebra 
\[
\frn_{5,2} = (0,0,0,12,13). 
\]
Here, the notation means that there exists a basis $(E^i)_{1\leq i\leq7}$ of $\frn_{5,2}^*$ such that $dE^4= E^1 \wedge E^2$, 
$dE^5= E^1 \wedge E^3$, and $dE^k=0$ otherwise. 
In particular, the Lie algebra of $\H$ is isomorphic to the one-dimensional extension $(\R\oplus\frn_{5,2})\rtimes \R$ of the six-dimensional 
2-step nilpotent Lie algebra $\R\oplus\frn_{5,2}$ by means of a certain derivation of $\R\oplus\frn_{5,2}$. Thus, $\H$ is almost nilpotent. 

\smallskip

From the computation of the canonical 1-form $h^{-1}dh$, we obtain the following basis of left-invariant 1-forms  on $\H$
\begin{equation}\label{Lcoframe}
\begin{split}
e^1 	&= \exp(x_7)\left(x_6 dx_7 +dx_6\right),\\
e^2	&= -\exp(x_7)\left(x_5 dx_7 +dx_5\right),\\
e^3	&= \exp(-x_7)\,dx_2,\\
e^4	&= \exp(x_7)\,dx_1,\\
e^5	&= 2\exp(x_7)\left(x_6\exp(x_7)\,dx_1+dx_4\right),\\
e^6	&= 2\exp(x_7)\left(x_5\exp(x_7)\,dx_1+dx_3\right),\\
e^7	&= -dx_7.
\end{split}
\end{equation}
We shall denote by $\mathcal{B}=(e_1,\ldots,e_7)$ the basis of left-invariant vector fields on $\H$ whose dual basis is $\mathcal{B}^*=(e^1,\ldots,e^7)$.
From \eqref{Lcoframe}, we see that
\begin{equation}\label{StrEqs}
\begin{split}
de^1&=0=de^2=de^7,\quad de^3=-e^{37},\quad de^4=e^{47},\\
de^5&= 2\,e^{14}+e^{57},\quad de^6 =  -2\,e^{24} + e^{67}. 
\end{split}
\end{equation}
Consequently, the non-zero Lie brackets between the left-invariant vector fields of $\mathcal{B}$ are the following
\begin{equation}\label{LieBr}
\begin{split}
&[e_3,e_7] = e_3,\quad [e_4,e_7] = -e_4,\quad [e_1,e_4]=-2e_5, \\
&[e_5,e_7] = -e_5,\quad [e_2,e_4] = 2e_6,\quad [e_6,e_7]=-e_6. 
\end{split}
\end{equation}

The left-invariant 3-form 
\[
\f = e^{127} + e^{347} + e^{567} + e^{135} - e^{146} - e^{236} - e^{245}
\]
defines a G$_2$-structure on $\H$ inducing the left-invariant metric $g_\f = \sum_{i=1}^7e^i\otimes e^i$ 
and the volume form $\mathrm{vol}_{g_\f}=e^{1234567}$. 
Using \eqref{StrEqs}, one can easily check that $\f$ is closed and that its intrinsic torsion form is
\[
\tau = -\star_\f d \star_\f \f = 2\,e^{12} + 2\,e^{34} -4\,e^{56}. 
\]
Moreover, the Hodge Laplacian of $\f$ is given by
\[
\Delta_\f\f = d\tau =  -8\left(e^{146}+e^{245}-e^{567}\right).
\] 
The closed G$_2$-structure $\f$ is a steady soliton of gradient type. 
Indeed, there exists a left-invariant gradient vector field $X$ on $\H$ for which  
\[
\Delta_\f\f = \mathcal{L}_{X}\f = d(\iota_{X}\f),
\]
namely 
\[
X = -4 e_7 = -4\left(e^7\right)^\sharp = \left(d(4x_7)\right)^\sharp = \nabla f, 
\]
where $f:\H\to\R$, $f(h) = 4x_7+b$, for some $b\in\R$. 

\smallskip

From the expression of $f$ and the Lie group isomorphism $\H\cong (\R\times \Ng_{5,2}) \rtimes \R$, 
it follows that this example is either an instance of case 1 or an instance of case 2(a) 
of the structure theorem for homogeneous gradient Laplacian solitons \cite[Theorem 1.1]{Ng}. 
To establish the right class, it is sufficient to compute the divergence of the symmetric 2-covariant tensor $t_\f\coloneqq g_\f(T_\f^2\cdot,\cdot)$, where 
$T_\f$ denotes the $g_\f$-skew-symmetric endomorphism corresponding to the torsion form $\tau$ via the identity $\tau = g_\f(T_\f\cdot,\cdot)$. 

Now, from the expressions of $\tau$ and $g_\f$, it follows that
\[
T_\f e_1=2e_2,~T_\f e_2 = -2e_1,~ T_\f e_3 = 2e_4,~ T_\f e_4 = -2e_3,~ T_\f e_5 = -4e_6,~ T_\f e_6 = 4e_5,~ T_\f e_7=0. 
\]
We compute the divergence of $t_\f$ using the formula 
\[
\mathrm{div}\left(t_\f\right)(\cdot) = \sum_{i=1}^7 g_\f\left(\nabla_{e_i}(T_\f^2(e_i)) -T_\f^2(\nabla_{e_i}e_i),\cdot\right). 
\]
Since $T_\f$ is left-invariant  and the matrix representation of $T_\f^2$ with respect to the basis $\mathcal{B}$ is diagonal, 
we only need to compute the covariant derivatives 
$\nabla_{e_k}e_k$, for $1\leq k \leq 7$. 
Since $(e_k)_{1\leq k\leq 7}$ is a left-invariant frame and the metric $g_\f$ is left-invariant, the Koszul formula gives 
\[
g_\f(\nabla_{e_i}e_j,e_k) = -\frac12\left(g_\f(e_i,[e_j,e_k]) + g_\f(e_j,[e_i,e_k]) +g_\f(e_k,[e_j,e_i]) \right),
\]
for all $1\leq i,j,k\leq 7$. Using this and \eqref{LieBr}, we obtain
\[
\begin{split}
\nabla_{e_k}e_k&=0,\quad k=1,2,7,\\
\nabla_{e_3}e_3&=-e_7,\\
\nabla_{e_k}e_k&=e_7,\quad k=4,5,6.
\end{split}
\]
Consequently, we have
\[
\begin{split}
\mathrm{div}\left(t_\f\right)(e_k)&=0,\quad 1\leq k \leq 6,\\
\mathrm{div}\left(t_\f\right)(e_7)&=-32.
\end{split}
\]
Since $\mathrm{div}\left(t_\f\right)\neq0$, the homogeneous space $(\H,g_\f)$ is a one-dimensional extension in the sense of \cite{PeWi}, and thus 
an instance of case 2(a) of \cite[Theorem 1.1]{Ng}. 

\bigskip

{\bf Acknowledgements.}
The authors were supported by GNSAGA of INdAM and by the project PRIN 2017  
``Real and Complex Manifolds: Topology, Geometry and Holomorphic Dynamics''. 
The first named author~is also supported by a grant from the Simons Foundation ($\#$944448). 
She would like to thank the Simons Foundation for the  support  to the  participation at  the  MATRIX Research Program
``Spectrum and Symmetry for Group Actions in Differential Geometry II'',  July 24 -- August 4, 2023.

\end{document}